\documentclass[a4paper]{amsart}

\usepackage{amsmath,amssymb,amsthm,amsfonts}

\newtheorem{thm}{Theorem}[section]
\newtheorem{lem}[thm]{Lemma}
\newtheorem{cor}[thm]{Corollary}
\newtheorem{const}[thm]{Construction}
\newtheorem{prop}[thm]{Proposition}
\newtheorem{deff}[thm]{Definition}
\newtheorem{rk}[thm]{Remark}
\newtheorem{examp}[thm]{Example}
\newtheorem{quest}[thm]{Question}
\newtheorem{claim}{Claim}
\newtheorem{prob}[thm]{Problem}

\newcommand{\brac}[1]{\langle #1\rangle}
\newcommand{\qedd}{\hfill \qedsymbol \vspace{1.5mm}}
\newcommand{\vs}{\vspace{0.3cm}}

\newcommand{\0}{\mathbf{0}}

\newcommand{\Z}{{\ensuremath{\mathbb{Z}}}}
\newcommand{\Fp}{\mathbb{F}_p}
\newcommand{\F}[1]{\mathbb{F}_{#1}}
\newcommand{\M}[1]{\mathcal{#1}}
\newcommand{\al}{\alpha}
\newcommand{\Ga}{\Gamma}
\newcommand{\AGaL}{{\rm A\Ga L}}
\newcommand{\GaL}{{\rm \Ga L}}
\newcommand{\AGL}{{\rm AGL}}
\newcommand{\GL}{{\rm GL}}
\newcommand{\SL}{{\rm SL}}
\newcommand{\PGL}{{\rm PGL}}
\newcommand{\PSL}{{\rm PSL}}
\newcommand{\Paley}{{\rm Paley}}
\newcommand{\GPaley}{{\rm GPaley}}
\newcommand{\oline}[1]{\overline{#1}}

\newcommand{\imply}{\Longrightarrow}
\newcommand{\mapp}{\longrightarrow}
\newcommand{\om}{\omega}
\newcommand{\omb}{\widehat{\omega}}

\newcommand{\Aut}{{\rm Aut}}
\newcommand{\Sym}{{\rm Sym}}
\newcommand{\Cay}{{\rm Cay}}
\newcommand{\Cyc}{{\rm Cyc}}

\newcommand{\GP}{\GPaley(p^R, \frac{p^R-1}{k})}
\newcommand{\GPq}{\GPaley(q, \frac{q-1}{k})}
\newcommand{\GPP}[2]{\GPaley(#1, \frac{#1 - 1}{#2})}

\newcommand{\bpf}{\begin{proof}}
\newcommand{\epf}{\end{proof}}

\newcommand{\bpro}{\begin{prob}}
\newcommand{\epro}{\end{prob}}

\newcommand{\bt}{\begin{thm}}
\newcommand{\et}{\end{thm}}

\newcommand{\bl}{\begin{lem}}
\newcommand{\el}{\end{lem}}

\newcommand{\bp}{\begin{prop}}
\newcommand{\ep}{\end{prop}}

\newcommand{\bc}{\begin{cor}}
\newcommand{\ec}{\end{cor}}

\newcommand{\bct}{\begin{const}}
\newcommand{\ect}{\end{const}}

\newcommand{\bdeff}{\begin{deff}}
\newcommand{\edeff}{\end{deff}}

\newcommand{\brk}{\begin{rk}}
\newcommand{\erk}{\end{rk}}

\newcommand{\bexp}{\begin{examp}}
\newcommand{\eexp}{\end{examp}}

\newcommand{\bq}{\begin{quest}}
\newcommand{\eq}{\end{quest}}

\newcommand{\bcl}{\begin{claim}}
\newcommand{\ecl}{\end{claim}}

\newcommand{\be}{\begin{eqnarray*}}
\newcommand{\ee}{\end{eqnarray*}}

\newcommand{\ben}{\begin{eqnarray}}
\newcommand{\een}{\end{eqnarray}}

\newcommand{\bi}{\begin{itemize}}
\newcommand{\ei}{\end{itemize}}

\newcommand{\bnum}{\begin{enumerate}}
\newcommand{\enum}{\end{enumerate}}

\begin{document}

\title[Generalised Paley Graphs]
{On generalised Paley graphs and\\ their automorphism groups}

\author[Tian Khoon Lim and Cheryl E. Praeger]
{Tian Khoon Lim\mbox{$^1$} and Cheryl E. Praeger}
\address{School of Mathematics and Statistics (M019)\\ The University of Western Australia\\
35 Stirling Highway\\ Crawley WA 6009\\ Australia} \email{limtk@phillip.com.sg, praeger@maths.uwa.edu.au}

\maketitle

\footnotetext[1]{Current address: 45 Choa Chua Kang Loop
\#07-14 Singapore 689679.}


\section{Introduction} \label{sec:intro}

The generalised Paley graphs are, as their name suggests, a generalisation of 
the Paley graphs, first defined by Paley in 1933 (see \cite{Paley}). 
They arise as the relation graphs of 
symmetric cyclotomic association schemes. However, their automorphism 
groups may be much larger than the groups of the corresponding schemes. 
We determine the parameters for which the graphs are connected, 
or equivalently, the schemes are primitive. Also we 
prove that generalised Paley graphs are sometimes isomorphic to 
Hamming graphs and
consequently have large automorphism groups, and we determine precisely the
parameters for this to occur. We prove that in the connected, non-Hamming case,
the automorphism group of a generalised Paley graph 
is a primitive group of affine type, and we find 
sufficient conditions under which the group is equal to the one-dimensional 
affine group of the associated cyclotomic association scheme. The results
have been applied in~\cite{LLP} to distinguish between cyclotomic schemes and similar 
twisted versions of these schemes, in the context of homogeneous factorisations
of complete graphs.

Let $\F{q}$ be a finite field with $q$ elements such that $q
\equiv 1$ (mod 4). Let $\om$ be a primitive element in $\F{q}$ and $S$ 
the set of nonzero squares in $\F{q}$, so $S=\{\om^2, \om^4, \ldots,
\om^{q-1}=1 \} =-S$. The \textit{Paley graph}, denoted by $\Paley(q)$, is 
the graph with vertex set $\F{q}$ and edges all pairs $\{x,y\}$ such that 
$x-y\in S$. The class of Paley graphs is one of the two infinite families 
of self-complementary arc-transitive graphs characterised by W. Peisert 
in \cite{Peisert2001}. Moreover, Paley graphs are also
examples of distance-transitive graphs, of strongly regular graphs, and of
conference graphs, see \cite[Section 10.3]{GR}. The automorphism group
$\Aut(\Paley(q))$ of $\Paley(q)$ is of index $2$ in the affine group
$\AGaL(1,q)$ and each permutation in $\AGaL(1,q) \setminus \Aut(\Paley(q))$
interchanges $\Paley(q)$ and its complementary graph. 
The generalised Paley graphs are defined similarly.

\bdeff \label{defn:gpaley} \rm{\textbf{(Generalised Paley Graph)} \ 
Let $\F{q}$ be a finite field of order $q$, and let $k$ be a divisor of 
$q-1$ such that $k\geq2$, and if $q$ is odd then $\frac{q-1}{k}$ is even. Let
$S$ be the subgroup of order $\frac{q-1}{k}$ of the multiplicative group 
$\F{q}^*$. Then the \textit{generalised Paley graph} $\GPq$ of $\F{q}$ 
is the graph with vertex set $\F{q}$ and edges all pairs 
$\{x,y\}$ such that $x-y\in S$. }\edeff

The generalised Paley graphs $\GPq$ are the relation graphs of the symmetric cyclotomic 
association scheme $\Cyc(q,k)$ defined in Subsection~\ref{sec:cyc}. They also
arise as the factors of cyclotomic homogeneous factorisations of complete 
graphs~\cite[Theorem 1.1]{LLP}. Our main theorem determines precise 
conditions under which  $\Cyc(q,k)$ is primitive (defined in 
Subsection~\ref{sec:cyc}), and shows how the automorphism group of a
generalised Paley graph depends heavily on the parameters $k$ and $q$.

\bt \label{thm:CYC} Let $V=\F{q}$, where $q=p^R$ with $p$ a prime, and let
$k\,\mid\,(q-1)$ such that $k>1$ and 
either $q$ is even or $\frac{q-1}{k}$ is even.
Let {\rm Cyc}$(q,k)$ be a $k$-class 
symmetric cyclotomic scheme with relation graphs   
$\Gamma_1, \ldots, \Gamma_k$, and let $\Ga=\GPq$. 
Then $\Ga_i\cong\Ga$ for each $i$ and: 
\bnum
\item[(1)] {\rm Cyc}$(q,k)$ is primitive if and only if 
$k$ is not a multiple of
$\frac{q-1}{p^a-1}$ for any proper divisor $a$ of $R$.
\item[(2)] $\Ga$ is a Hamming graph if and only if $k=\frac{a(q-1)}
{R(p^a-1)}$ for some proper divisor $a$ of $R$.
\item[(3)] If $\Ga$ is connected and is not a Hamming graph, then  
$\Aut(\Gamma)$ is a primitive subgroup of $\AGL(R,p)$ containing the 
group of translations $\Z_p^R$.
\item[(4)] If $k$ divides $p-1$, then $\Aut(\Gamma)=
\Aut({\rm Cyc}(q,k))< \AGaL(1,q)$. \enum \et

Commentary on the significance and consequences of this result 
is given in Subsection~\ref{sub:aut}. In particular, Hamming graphs 
and their automorphism groups are discussed there. The theorem may
be contrasted with McConnel's Theorem~\cite{MC}, proved in 1963, 
that the automorphism group of $\Cyc(q,k)$ (the intersection of 
the automorphism groups of its relation graphs) is always a subgroup 
of $\AGaL(1,q)$.

\subsection{Cyclotomic association schemes and cyclotomic factorisations} 
\label{sec:cyc}
Here we describe briefly 
the relationship between generalised Paley 
graphs, symmetric cyclotomic association schemes, and cyclotomic homogeneous factorisations of complete graphs. Symmetric association 
schemes are defined as follows.

\bdeff \label{def:AS} \cite[p. 43]{bcn} \rm{A \textit{symmetric $k$-class 
association scheme} is a
pair $(V,\M{R})$ such that \bnum
\item $\M{R} =\{R_0, R_1, \ldots, R_k\}$ is a partition of $V \times V$;
\item $R_0=\{(x,x) \mid x\in V\}$;
\item $R_i=R_i^T$ (that is, $(x,y) \in R_i$ $\imply$ $(y,x) \in R_i$) for all 
$i \in\{0,1,\ldots, k\}$;
\item there are constants $p_{ij}^h$ (called the \textit{intersection numbers} 
of the scheme) such that for any pair $(x,y) \in R_h$, the number of 
elements $z \in V$ with $(x,z) \in R_i$ and
$(z,y) \in R_j$ equals $p_{ij}^h$. \enum} \edeff

For each $i\geq1$, the class $R_i$ corresponds to the undirected
graph $\Gamma_i =(V, E_i)$, where $E_i =\{\{x,y\} \mid (x,y) \in
R_i\}$ (see for example \cite[Chapter 12]{Godsil93}). These graphs are called
the \emph{relation graphs} of the scheme and are not in general isomorphic. 
The scheme is said to be 
\textit{primitive} if each of the $\Gamma_i$ is connected, and
otherwise it is called \textit{imprimitive}. The automorphism group 
$\Aut(V,\M{R})$ is the largest subgroup of $\Sym(V)$ that, in its 
natural action on $V\times V$, fixes each of the relations $R_1,\dots,R_k$
setwise, that is to say,  $\Aut(V,\M{R})=\bigcap_{i=1}^k\Aut(\Ga_i)$.
The edge sets of the $\Ga_i$ 
form a partition $\M{E}=\{E_1,\dots,E_k\}$ of the edge set of the 
complete graph $K_n$, where $|V|=n$, and hence the relation graphs form a
factorisation of $K_n$.  If the subgroup of $\Sym(V)$ fixing $\M{R}$ 
setwise permutes transitively the set $\{\Ga_1,\dots,\Ga_k\}$ and if
$\Aut(V,\M{R})$ is transitive on $V$, 
then the factorisation is called \emph{homogeneous}. In particular, in
this case the relation graphs $\Gamma_i$ are pairwise isomorphic.

Let $V=\F{q}$, let $k \mid (q-1)$ such that $k>1$ and  either $q$ is even or 
$\frac{q-1}{k}$ is even, and let $S(k) =\brac{\om^k} \subseteq V^* = 
\F{q}^*$ (the multiplicative group of $\F{q}$).
Then the \emph{$k$-class symmetric cyclotomic scheme} 
Cyc$(q,k)=(V, \M{R})$, has $R_i=\{(x,y) \mid y-x \in S(k)\om^i\}$ 
for $1 \leq i \leq k$. Note that condition (3) of Definition~\ref{def:AS} holds
since $-1\in S(k)$.  Also, the relation graph $\Ga_k$ is the generalised Paley 
graph of $\F{q}$ relative to $S(k)$, and in particular if $k=2$, then $\Ga_2$ 
is the Paley graph of $\F{q}$ (see \cite[p. 66]{bcn}). Moreover, the affine 
group $\AGL(1,q)$ fixes $\M{R}$ setwise and permutes the relation
graphs transitively, and $\Aut(\Cyc(q,k))=\cap_{i-1}^k\Aut(\Ga_i)$ contains 
the group of translations which is transitive on $V$. 
Thus the relation graphs for
Cyc$(q,k)$ form a homogeneous factorisation of $K_q$ 
called a \emph{cyclotomic factorisation}. In particular all of the relation 
graphs are isomorphic to $\GPq$, and $\Aut(\GPq)$ contains the subgroup of
$\AGL(1,q)$ of order $q(q-1)/k$ acting arc-transitively.

\subsection{Commentary on Theorem~\ref{thm:CYC}: automorphism groups of graphs and schemes}
\label{sub:aut}
In 1963, McConnel~\cite{MC} proved that the full automorphism 
group of $\Cyc(q,k)$ is a subgroup of $\AGaL(1,q)$. 
However, the automorphism groups of the relation graphs
of $\Cyc(q,k)$, that is, of the generalised Paley graphs $\GPq$, 
may be much larger. This paper initiates
a study of these automorphism groups for various ranges of values of 
the parameters $q$ and $k$. 

The `easiest' way for the automorphism group of $\GPq$ to be larger is if
the graph is not connected. Since all the relation graphs of $\Cyc(q,k)$ 
are isomorphic to $\GPq$, it follows that $\Cyc(q,k)$ is primitive if
and only if $\GPq$ is connected. We determine in Theorem~\ref{thm:CYC}~(1) 
the precise parameter values for $\GPq$ to be connected, or equivalently, 
for $\Cyc(q,k)$ to be primitive. Moreover in 
Theorem~\ref{thm:add}, a more detailed version of this result, we prove
that, if $\GPq$ is disconnected, then 
its connected components are generalised Paley graphs over a 
proper subfield and generate a cyclotomic scheme over this subfield.

Suppose now that $\GPq$ is connected. It is possible for $\GPq$ to be a 
Hamming graph, for certain $q$ and $k$, and hence 
have a large automorphism group. For positive
integers $a>1, b>1$, the \textit{Hamming graph} $H(a,b)$ has as vertices all 
$b$-tuples with entries from a set $\Delta$ of size $a$. Two vertices are 
adjacent in $H(a,b)$ if and only if the two $b$-tuples differ in exactly 
one entry. The full automorphism group of $H(a,b)$ is the wreath product 
$S_a \wr  S_b$ in its product action on $\Delta^b$, see 
\cite[Theorem 9.2.1]{bcn}. Section~\ref{sec:prelim} contains details about 
the product action. We determine in Theorem~\ref{thm:CYC}~(2) 
the precise parameter values for $\GPq$ to be a Hamming graph. 

Moreover we 
prove, in Theorem~\ref{thm:CYC}~(3), that if $\GPq$ is connected and not a 
Hamming graph, then its automorphism group is a primitive group of affine 
type. Although this gives a lot of information about the group, it does not 
determine it completely. In a special case, that is relevant to our work
on homogeneous factorisations in \cite{LLP}, we were able to show that
the automorphism group of $\GPq$ is indeed equal to 
the automorphism group of Cyc$(q,k)$ (and no larger), 
see Theorem~\ref{thm:CYC}~(4).
 
We make additional detailed comments about Theorem~\ref{thm:CYC} and 
its consequences in Remark~\ref{rk:appn}.

\begin{rk}\label{rk:appn}{\rm
(a) The graph $\GPq$ is a Cayley graph for the translation subgroup $T$
of $\AGaL(1,q)$, see Subsection~\ref{sub:gpcay}. Theorem~\ref{thm:CYC}~(3) 
proves that, provided $\GPq$ is connected and not a Hamming graph, then
it is a \textit{normal Cayley graph}, that is,
the translation subgroup $T$ of automorphisms is normal in the full 
automorphism group. 

(b)   The condition $k \mid (p-1)$ holds in particular 
if $q=p$, that is to say, if $R=1$. In this case, Theorem~\ref{thm:CYC}~(4) 
follows from an old result of W. Burnside about primitive 
permutation groups of prime degree, see \cite[11.7]{wie}. 

(c) The proof of Theorem~\ref{thm:CYC}~(3) uses results from \cite{Praeger90,
PS92}, and that of Theorem~\ref{thm:CYC}~(4) depends heavily on results 
in \cite{GPPS}. As the results used from these papers rely on the classification
of the simple groups, these two parts of Theorem~\ref{thm:CYC} also
rely on that classification.

(d) Apart from the possibilities that $\GPq$ may be disconnected or
isomorphic to a Hamming graph, which are 
dealt with in Theorem~\ref{thm:CYC}~(1) and (2), there are other cases 
where $\Aut(\GPq)$ is not a one-dimensional affine group. 
In Example~\ref{ex:1}, we give an explicit example. Thus despite our
identifying the disconnected and Hamming cases precisely in
Theorem~\ref{thm:CYC}, there remain some mysteries to be solved concerning
generalised Paley graphs:  
Problem~\ref{prob:p1} below is still largely open. 

(e) Our interest in generalised Paley graphs arose from our study of 
homogeneous factorisations of complete graphs (see Subsection~\ref{sec:cyc}). 
These factorisations were introduced in \cite{LP2003} as a 
generalisation of vertex-transitive self-complementary graphs.
Our study in \cite{LLP} gave a classification of arc-transitive homogeneous 
factorisations of complete graphs. In addition to the cyclotomic
factorisations, we discovered a new family of examples that 
generalise an infinite family of vertex-transitive self-complementary 
graphs constructed and characterised by Peisert~\cite{Peisert2001}.
They may be viewed as a twisted version of cyclotomic factorisations. Using 
Theorem~\ref{thm:CYC}~(4), we proved in \cite{LLP} that factor graphs in 
these two homogeneous factorisations, with the same parameters $q$ and $k$, 
were non-isomorphic; we showed that their automorphism groups had 
non-isomorphic intersections with $\AGaL(1,q)$.
} 
\end{rk}

\begin{prob}\label{prob:p1} \rm{Determine the precise conditions on $k$ and $p^R$ under which the conclusion of Theorem~\ref{thm:CYC}~(4)
holds.}
\end{prob}

\bexp \label{ex:1} \rm{Take $R=4$, $p=3$ and $k =4$, so that $\frac{p^R-1}{k} 
=20 $, and let $\Gamma =\GPaley(81, 20)$. Then $k \neq
\frac{a(p^R-1)}{R(p^a-1)}$, and $k$ is not a multiple of $\frac{q-1}{p^a-1}$,
for any proper divisor $a$ of $R$. 
Hence, by Theorem~\ref{thm:CYC},
$\Gamma$ is connected and not a Hamming graph, and $\Aut(\Gamma)$ 
is a primitive subgroup of $\AGL(4,3)$. Using \textsc{Magma} \cite{Magma}, we
computed $\Aut(\Gamma)$. Its order is $|\Aut(\Gamma)| = 233 280$, greater 
than $|\AGaL(1,81)|=25 920$. Thus $\Aut(\Gamma)$ is not contained in the 
one-dimensional affine group. A further
check using \textsc{Magma} showed that a point-stabiliser $A_0$ of $A:=\Aut
(\Gamma)$, which has order $2880$, contains a normal subgroup $B$ 
isomorphic to $A_6$, and $A_0/B \cong D_8$, so $\Aut(\Gamma)
= Z_3^4 \rtimes ( A_6 \cdot D_8)$.}\eexp

In Section~\ref{sec:prelim}, we 
introduce some terminology and definitions needed for subsequent results, and
we prove Theorem~\ref{thm:CYC}~(1). We  prove Theorem~\ref{thm:CYC}~(2), 
(3) and (4) in Sections~\ref{sec:A}, \ref{sec:pfthm1} and~\ref{sec:subcase1}  
respectively.

\section{Preliminaries and Proof of Theorem~\ref{thm:CYC}~(1)}
\label{sec:prelim}

\subsection{Cayley graphs.}\label{sub:cay}
All graphs considered are finite, undirected and without loops or multiple 
edges. 
Thus a graph $\Gamma=(V,E)$ consists of a vertex set $V$ and a subset $E$ of 
unordered pairs from $V$, called the edge set. 
An arc is an ordered pair $(u,v)$ where $\{u,v\}$ is an edge.
The generalised Paley graphs belong to a larger class of graphs called the
Cayley graphs, defined as follows.

\bdeff \label{defn:cay} \rm{\textbf{(Cayley Graph)} \ For a group $K$ and a 
nonempty subset $H$ of $K$ such that $1_K \notin H$ and $H=H^{-1}=\{h^{-1} 
\mid h \in H\}$, the \textit{Cayley graph} $\Gamma = \Cay(K,H)$ of $K$ 
relative to $H$ is the graph with vertex set $K$ such that
$\{x,y\}$ is an edge if and only if $xy^{-1} \in H$.
} \edeff

A Cayley graph $\Gamma = \Cay(K,H)$ is connected if and only if $\brac{H} = K$.
Furthermore, if $\Gamma =\Cay(K,H)$ is disconnected, then each connected 
component is isomorphic to $\Cay(\brac{H},H)$ and the number of connected 
components equals $\frac{|K|}{|\brac{H}|}$.
A permutation group $G$ on $V$ is \textit{semiregular} if the only element 
fixing a point in $V$ is the identity element of $G$; and 
$G$ is \textit{regular} on $V$ if it is both semiregular
and transitive. In Definition~\ref{defn:cay}, the group $K$ acts regularly on 
vertices by $y:x\rightarrow xy$, for $x,y\in K$. Conversely, see
\cite[Lemma 16.3]{Biggs_b93}, a graph $\Gamma$ is isomorphic to a Cayley
graph for some group if and only if $\Aut(\Gamma)$ has a subgroup which is
regular on vertices. 

\subsection{Generalised Paley graphs as Cayley graphs.}\label{sub:gpcay}
In follows from Definition~\ref{defn:gpaley} that 
\[
\GPq=\Cay(V,S),
\] 
where $V$ is the additive group of the field $\F{q}$ and $S$ is the 
unique subgroup of order $\frac{q-1}{k}$ 
of the multiplicative group $\F{q}^*$. Let
$\om$ be a primitive element of $\F{q}$. Then $S=\brac{\om^k}$. Thus $\GPq$
admits the additive group of $\F{q}$ acting regularly by $t_y:x\mapp x+y$
(for $x,y\in\F{q}$)  as a subgroup of automorphisms. To distinguish this 
subgroup from the vertex set $V$, we denote it by $T=\{t_y | y\in\F{q}\}$, and 
call it the \emph{translation group} of $\AGaL(1,q)$. 

Now  $T\cong\Z_p^R$, where
$q=p^R$ with $p$ prime, and $T$ is the unique minimal normal subgroup of
$\AGaL(1,q)$.
Let $\omb$ denote the scalar multiplication map $\omb: x \mapp x\om$ (for all 
$x \in \F{q}$) corresponding to the primitive element $\om$, and  
let $\al$ denote the Frobenius automorphism of $\F{q}$,
that is, $\al: x \mapp x^p$. Then $\AGaL(1,q)=T\rtimes\brac{\omb,\al}$, and 
the one-dimensional general semilinear group $\GaL(1,q) =\brac{\omb, \al}$.
Now both $\omb^k$ and $\al$ fix $\0$ and fix $S$ setwise, and hence
$T \rtimes \brac{W,\al}$ is a subgroup of automorphisms of $\GPq$,
where $W:=\brac{\omb^k}$. 
In fact, $T \rtimes W$ is arc-transitive, and $W$ is transitive on 
$S = \{ 1, \om^k, \om^{2k}, \ldots, \om^{((q-1)/k )-k}\} =1^W$. 
(For a permutation group $K$
on $V$ and a point $v\in V$ we denote by $v^K$ the $K$-orbit $\{v^x\,|\,x\in 
K\}$ containing $v$.)

\subsection{Hamming graphs and Cayley graphs.}\label{sub:hamcay}
Let $H$ be a group, $b$ a positive integer and $K$ be a subgroup of the 
symmetric group $S_b$.
Then the \emph{wreath product} $H \wr  K$ is the semidirect product $H^b 
\rtimes K$ where elements of
$K$ act on $H^b$ by permuting the ``entries" of elements of $H^b$, that is, 
$(h_1, h_2, \ldots,h_b)^{k^{-1}} = (h_{1^k}, h_{2^k}, \ldots, h_{b^k})$ for 
all $(h_1, h_2, \ldots, h_b) \in H^b$ and
$k \in K$. Now suppose $H \leq$ Sym$(\Delta)$. Then the \textit{product 
action} of $H \wr \ K$ on
$\Delta^b$ is defined as follows. Elements of $H^b$ act coordinate-wise on 
$\Delta^b$ and elements
of $K$ permute the coordinates: for $(h_1, \ldots, h_b) \in H^b$, $k \in K$, 
and $(\delta_1,
\ldots, \delta_b) \in \Delta^b$, \be (\delta_1, \ldots, \delta_b)^{(h_1, 
\ldots, h_b)} & = &(\delta_1^{h_1}, \ldots, \delta_b^{h_b}) \\
(\delta_1, \ldots, \delta_b)^{k^{-1}} & = & (\delta_{1^k}, \ldots, \delta_{b^k}). \ee

If $A$ is a regular subgroup of $S_a$  
then $A^b < S_a \wr \ S_b$ and $A^b$ acts regularly on the vertices of
$H(a,b)$. Thus (see Subsection~\ref{sub:cay}) $H(a,b)$ is a Cayley graph. 
If $a$ is a prime power $q$, then $A$ can be identified with the additive 
group of a finite field $\F{q}$, and the vertex set of $H(a,b)$
can be identified with $\F{q}^b$.

\subsection{Primitive permutation groups}\label{sub:prim}
Let $G$ be a transitive permutation group acting on a finite set $V$. 
A nonempty subset $\Delta
\subseteq V$ is called a \textit{block} for $G$ if for every $g \in G$, 
either $\Delta \cap
\Delta^g = \emptyset$ or $\Delta = \Delta^g$. A block $\Delta$ is said to 
be \textit{trivial} if
$|\Delta| = 1$ or $\Delta = V$. Otherwise, $\Delta$ is called 
\textit{nontrivial}. We say that the group $G$ is \textit{primitive} if 
the only blocks for $G$ are the trivial ones.

The possible structures of finite primitive permutation groups up to 
permutational isomorphism are described by the O'Nan-Scott Theorem (for 
example see \cite{Cameron_b99} or \cite{LPS88}). Here, we will briefly 
describe the three types of finite
primitive permutation groups relevant to this paper (we refer readers to 
\cite{Cameron_b99, LPS88}
for further details about the remaining types).

A finite primitive permutation group $G$ on $V$ is of type {\sc HA} (holomorph 
of an abelian group) if $G = T \rtimes G_0$ is a subgroup of an affine group 
$\AGL(R,p)$ on $V$, where $T \cong\Z_p^R$ is the (regular) group of 
translations (and we may identify $V$ with $\Z_p^R$) and $G_0$
is an irreducible subgroup of $\GL(R,p)$. We often say that a primitive group 
of this type is of affine type.
A primitive permutation group $G$ is of type {\sc AS} if $G$ is an 
\emph{almost simple group}, that is $N \leq G \leq \Aut(N)$ where $N$ is a 
finite nonabelian simple group. Such a group can equivalently be defined as a 
primitive group $G$ having a unique minimal normal subgroup $N$ which is
nonabelian and simple.
Finally, a primitive permutation group $G$ on $V$ is of type {\sc PA} 
(product action) if $V=\Delta^b$ and $N^b \leq G \leq H \wr \ S_b \leq$ 
Sym$(\Delta) \wr \ S_b$ in its product action, where $H$ is a primitive 
permutation group on $\Delta$ of type {\sc AS} with simple
normal subgroup $N$.

\subsection{Proof of Theorem~\ref{thm:CYC}~(1)}
The following result Theorem~\ref{thm:add} relates the connectedness of 
$\Gamma$ with the action of $W$ on $V$, and
Theorem~\ref{thm:CYC}~(1) follows immediately from it.

\bt \label{thm:add} Let $\Gamma =\GPq =\Cay(V,S)$, where
$V =\F{q}$ and $S =\brac{\om^k}$, with $k$ a divisor of $q-1$ such that $k\geq
2$ and either $q$ or $\frac{q-1}{k}$ is even. Let $q=p^R$ with $p$ prime, and
let $\Cyc(q,k)=(V,\M{R})$.  
\bnum
\item The following are equivalent:
\begin{enumerate}
\item[(i)] $\Gamma$ is connected,
\item[(ii)] $\Cyc(q,k)$ is primitive,
\item[(iii)] $\brac{\omb^k}$ acts irreducibly on $V$, 
\item[(iv)] $k$ is not a multiple of $\frac{q-1}{p^a-1}$ for any proper 
divisor $a$ of $R$.
\end{enumerate}
\item Suppose that $k$ is a multiple of $\frac{q-1}{p^a-1}$, where $a$
divides $R$, so that $\Gamma$ is not connected. Then the connected components 
of $\Gamma$ are all isomorphic, and the component  $\Gamma_0$ containing $\0$ 
has vertex set $\F{p^a}$ (a proper subfield  of $\F{q}$)
containing $S$, and is isomorphic to $\GPaley(p^a, \frac{p^a-1}{k'})$, 
where $k' =\frac{p^a-1}{q-1} \ k\geq 1$.
Furthermore, $\Aut(\Gamma) =
\Aut(\Gamma_0) \wr \ S_{p^{R-a}}$. \enum \et 

We note that in part (2), $k'$ may equal $1$, and if this happens we 
still use the notation $\GPaley(p^a, \frac{p^a-1}{k'})$ for $\Ga_0$ even though
in this case $\Ga_0=\Cay(\F{p^a},S)\cong K_{p^a}$, the complete graph on 
$p^a$ vertices.

\bpf (1) \ The equivalence of parts (i) and (ii) follows from our discussion in
Subsection~\ref{sec:cyc}. Let $U$ be the $\Fp$-span of $S$, that is,
$U = \{ \sum_{\om^{ik} \in S} \lambda_i \om^{ik} \mid \lambda_i \in \Fp \}$. 
Since $W=\brac{\omb^k}$ leaves $S$ invariant, it also leaves invariant the 
$\Fp$-span $U$ of $S$. Also we note that $\Gamma$ is connected if and only if 
$U = V$ (see Definition~\ref{defn:cay}).

Suppose $W$ acts irreducibly on $V$. Then as $U$ is $W$-invariant and nonzero, 
$U = V$ and hence $\Gamma$ is connected. Conversely suppose $\Gamma$ is 
connected. Then $S$ is an $\F{p}$-spanning set for $V$ (that is $U=V$). 
Now $S$ is the orbit $1^W = \brac{\om^k}$ (note that $W$ acts by field
multiplication). Also, for each $\om^i \in V^*$, $\omb^i$ maps $S$ to 
$S \om^i$, and as $\omb^i\in \GL(1,p^R)$, it follows that $S \om^i$ is also an 
$\F{p}$-spanning set for $V$. 
However $S \om^i$ is the $W$-orbit containing $\om^i$. Hence 
every $W$-orbit in $V^*$ is a spanning set for $V$, and so $W$ is 
irreducible on $V$. Thus (i) and (iii) are
equivalent.

Suppose that $k$ is a multiple of $\frac{q-1}{p^a-1}$ for some proper 
divisor $a$ of $R$. Then $S$ is a subgroup of the multiplicative group
of the proper subfield $\F{p^a}$ of $\F{q}$. Thus $U\subset\F{p^a}$ and so
$\Ga$ is disconnected. Therefore condition (iv) implies condition (i). The
reverse implication will follow from (2).

(2) \ Suppose $\Gamma$ is disconnected. 
Let $U$ be the vertex set of the connected component of
$\Gamma$ containing $1 \in \F{p^R}$. Then $U$ is the $\Fp$-span of $S$. 
It follows that all the
connected components of $\Gamma$ are isomorphic to $\Cay(U,S)$. 
We claim that $U$ is a subfield of $V =\F{p^R}$.

Since $U$ is $W$-invariant, $U^{\omb^{ik}} =U$ for each $\omb^{ik} \in W$, 
and hence $U \om^{ik} =U$ for each $\om^{ik} \in S$. 
Thus $U$ is closed under multiplication by elements of $S$. Now each
element of $U$ is of the form $\sum_{\om^{ik} \in S} \lambda_i \om^{ik}$ 
for some $\lambda_i \in\Fp$, and (by regarding $\lambda_i$ as an integer 
in the range $0 \leq \lambda_i \leq p-1$) each
$\lambda_i \om^{ik}$ is equal to the sum $\om^{ik} + \cdots + \om^{ik}$ 
($\lambda_i$-times). Thus
each element of $U$ is a sum of a finite number of elements of $S$. 
Since $U$ is closed under
addition and under multiplication by elements of $S$, it follows that 
$U$ is closed under multiplication. 
Thus $U$ is a subring of $V$. 
Also $U$ contains the identity $1$ of $V$ (since $1\in S$). 
Let $u \in U \setminus \{0\}$. Then since $V$ is finite, $u^i =u^{i+j}$ for
some $i \geq 1$ and $j \geq 1$, and hence $u^{-1} = u^{j-1} \in U$. 
Thus $U$ is a subfield of $V=\F{p^R}$ as claimed.

Hence $|U| =p^a$ for some proper divisor $a$ of $R$. Also, since $S
\leq U^*$, it follows that $|S| = \frac{p^R-1}{k}$ divides $p^a-1$. Let $k' =
\frac{(p^a-1)k}{p^R-1}$. Then $|S| = \frac{p^R-1}{k} =\frac{p^a-1}{k'}$, 
and by definition of a generalised Paley graph, we have $\Cay(U,S) = 
\GPaley(p^a, \frac{p^a-1}{k'})$ (though perhaps $k'=1$). Since there are
$p^{R-a}$ connected components in $\Gamma$, it follows that $\Aut(\Gamma) = 
\Aut(\GPaley(p^a,\frac{p^a-1}{k'})) \wr \ S_{p^{R-a}}$. \epf

From Theorem~\ref{thm:add}, if $\Gamma =\GPq$
is disconnected, then the connected components are generalised Paley graphs 
for subfields. In the rest of the paper we will assume that $\Gamma =\GPq$ 
is connected.

\section{Proof of Theorem~\ref{thm:CYC}~(4)} \label{sec:A}

Suppose first that $\Gamma = \GP \cong H(p^a,b)$, where $R =ab$ with $b>1$. 
The valency of $\Gamma$ is
$\frac{p^R -1}{k} = b(p^a-1)$, so $k = \frac{p^R-1}{b(p^a-1)} =
\frac{a(p^R-1)}{R(p^a-1)}$ as required.
Conversely suppose that $k = \frac{a(p^R-1)}{R(p^a-1)}$ where $R=ab$ and 
$1 \leq a < R$. Then since $\Gamma = \GP$ is connected, the $\Fp$-span of 
$S =\brac{\om^k}$ equals $V$, that is, $\{ \sum_{\om^{ik} \in S} \lambda_i 
\om^{ik} \mid \lambda_i \in \Fp \} =V$. Let $U$
be the $\F{p^a}$-span of the set $X := \{1, \om^k, \om^{2k}, \ldots, 
\om^{(b-1)k} \}$. We claim that $U=V$.

Now $\F{p^a}, X \subseteq V =\F{p^R}$, and hence $U \subseteq V$. 
Moreover, since $k =
\frac{a(p^R-1)}{R(p^a-1)}$, the set $S =\brac{\om^k}$ has order 
$\frac{p^R-1}{k} = \frac{R}{a}
\cdot (p^a-1) = b(p^a-1)$, and also $\om^{bk}$ has order $p^a-1$.  
Thus $\brac{\om^{bk}} =
\F{p^a}^*$. Now suppose $v \neq 0$ and $v \in V$. Then $v = \sum 
\lambda_i \om^{ik}$ where
$\lambda_i \in \Fp$ (not all zero) and the sum is over all $\om^{ik} \in S$. 
Let $i =bx_i +r_i$ where $0 \leq r_i < b$. 
Then $\om^{ik} = \om^{(bx_i +r_i)k} = \om^{bx_ik} \cdot \om^{r_ik}$ and we
have 
\be 
v \ = \ \sum \lambda_i \om^{ik} & = & \sum \lambda_i \om^{bx_ik} 
\cdot \om^{r_ik}. 
\ee
Since $\F{p^a}^* =\brac{\om^{bk}}$, it follows that $\om^{bx_ik} \in \F{p^a}$. 
Thus $\lambda_i
\om^{bx_ik} \in \F{p^a}$ and so $v \in U$. Thus $U=V$ as claimed.

From now on we shall regard $V$ as a vector space over $\F{p^a}$. 
We have shown that $V$ is
spanned by the set $X = \{1, \om^k, \om^{2k}, \ldots, \om^{(b-1)k} \}$, 
and as dim$_{\F{p^a}}(\F{p^R}) =b$, it follows that $X$ is an $\F{p^a}$-basis 
for $V$. 
Define $\Theta: V\mapp \F{p^a}^b$ as follows. For $u  = \sum_{j=0}^{b-1} 
\mu_j \om^{jk}$ with $\mu_j \in \F{p^a}$,
let $ \Theta(u) = (\mu_0, \mu_1, \ldots, \mu_{b-1}) \in \F{p^a}^b$. 
Since $X$ is an
$\F{p^a}$-basis for $V$, $\Theta$ is a bijection.

Next, we determine the image of the connecting set $S \subset V$ under $\Theta$. As we observed
above, $|S| = |\brac{\om^k}| =b(p^a-1)$. Thus each element of $S$ can be expressed uniquely as
$\om^{ik}$ for some $i$ such that $0 \leq i \leq b(p^a-1)-1$. As before, we write $i =bx_i +r_i$
where $0 \leq r_i \leq b-1$. Thus $\om^{ik} =\om^{bx_ik} \cdot \om^{r_ik}$. Now $\om^{bx_ik} \in
\F{p^a}^* =\brac{\om^{bk}}$, and so 
\[ 
\Theta(\om^{ik}) = \Theta(\om^{bx_ik} \cdot \om^{r_ik}) =
(0, \ldots, 0, \underbrace{\om^{bx_ik}}_{r_i{\rm th}}, 0, \ldots, 0). 
\] 
Observe that $x_i$ can be
any integer satisfying $0 \leq x_i \leq p^a-2$, and therefore $\om^{bk x_i}$ takes on each of the
values in $\F{p^a}^*$. Moreover each of these values occurs exactly once in each of the positions
$r$, for $0 \leq r \leq b-1$. Thus $\Theta (S)$ is the set of all elements of $\F{p^a}^b$ with
exactly one component non-zero, that is, the set of ``weight-one" vectors.

Now $\Theta$ determines an isomorphism from $\Gamma$ to the Cayley graph for $\F{p^a}^b$ with
connecting set $\Theta(S)$. In this Cayley graph, two $b$-tuples $u, v \in \F{p^a}^b$ are adjacent
if and only if $u-v \in \Theta(S)$, that is, if and only if $u-v$ has exactly one non-zero
component. Thus $\Gamma =\GP$ is mapped under the isomorphism $\Theta$ to the Hamming graph
$H(p^a, b)$ where $b =\frac{R}{a}$. \qedd

\section{Proof of Theorem~\ref{thm:CYC}~(3)} \label{sec:pfthm1}

Recall that by Theorem~\ref{thm:add}~(1), if $\Gamma=\GPq$ is connected, 
then $W=\brac{\omb^k}$ acts irreducibly on $V$. 
It follows that the group $G = T \rtimes W$ is a vertex-primitive
subgroup of $\Aut(\Gamma)$ of affine type. Thus $\Aut(\Gamma)$ is a primitive
permutation group on $V$ containing $G$. 
We will use results from \cite{Praeger90,PS92} concerning such groups.

\vs \noindent \textit{Proof of Theorem~\ref{thm:CYC}~(3).} \ Suppose that
$\Gamma=\GP =\Cay (V,S)$ is connected and not a Hamming graph. 
Let $G = T \rtimes W $ where $T \cong \Z_p^R$ and 
$W =\brac{\omb^k}$, as in Subsection~\ref{sub:gpcay}. By 
Theorem~\ref{thm:CYC}~(2), $k \neq \frac{a(p^R-1)}{R(p^a-1)}$ for any $a 
\mid R$ with $1 \leq a < R$.
Also, by Subsection~\ref{sub:gpcay} and Theorem~\ref{thm:add}, 
$G \leq X:=\Aut(\Gamma)$, and $W$ is irreducible on $V$, so  
$G$ is a primitive subgroup of $\AGL(R,p)$.
Suppose, for a contradiction, that $X$ is not contained in $\AGL(R,p)$.
Since $k\geq 2$, $\Gamma$ is not a complete graph and so $X \neq S_{p^R}$ or 
$A_{p^R}$. By \cite[Proposition 5.1]{Praeger90}, it follows that $X$ is 
primitive of type {\sc PA}. Thus (see
Subsection~\ref{sub:prim}) $R=ab$ with $b \geq 2$, $V = \Delta^b$ where 
$|\Delta|=p^a$, and $N^b\leq X \leq H \wr \ S_b$ with $H$ primitive on 
$\Delta$ of type {\sc AS} with simple normal subgroup $N$. 
Moreover from \cite[Proposition 5.1]{Praeger90} and
\cite[Proposition 2.1]{PS92}, either $N = A_{p^a}$ or $N$ and $p^a$ are as 
listed in \cite[Table 2]{Praeger90} (denoted as $L_1$ in \cite{Praeger90}). 
In all cases $N$ acts 2-transitively on $\Delta$, and since $N$ is nonabelian 
simple, it follows that $|\Delta|=p^a\geq5$.

We will prove that $\frac{p^R-1}{k} = b(p^{a}-1)$, contradicting the 
assumption on $k$, and therefore proving the theorem. Let $\gamma \in \Delta$ 
and consider the point $u=(\gamma, \ldots, \gamma) \in
\Delta^b = V$. Since $N^b$ is transitive on $V$, $X=N^b X_{u}$, where $X_u$ is 
the stabiliser of $u$ in $X$. Now $X_{u}$ contains $(N^b)_{u} = 
(N_{\gamma})^b$, and $N_{\gamma}$ is
transitive on $\Delta - \{\gamma\}$. If $v \neq u$, then 
$v := (\delta_1, \ldots, \delta_b)$ with, say, $\ell$ entries different from 
$\gamma$, where $1 \leq \ell \leq b$, and the length of the $(N^b)_{u}$-orbit 
containing $v$ is $(p^{a}-1)^{\ell}$.

Since $X$ is a primitive subgroup of Sym$(\Delta) \wr \ S_b$, $X$ projects to 
a transitive subgroup of $S_b$ (for instance, see 
\cite[Theorem 4.5]{Cameron_b99}). Moreover, since $X=N^bX_u$, it follows that
$X_u$ also projects to a transitive subgroup of $S_b$. Thus 
the $\ell$-subset of subscripts $i$ such that $\delta_i
\neq \gamma$ has $n_{\ell}$ distinct images under $X_u$, where 
$n_{\ell} \geq b/\ell$. It follows that the length of the $X_{u}$-orbit 
containing $v$ is at least $n_{\ell} \cdot
(p^{a}-1)^{\ell} \geq \frac{b}{\ell} \cdot (p^{a}-1)^{\ell}$. Suppose now that the point $v$ has
been chosen to lie in $\Gamma(u)$ (the set of all vertices in $\Gamma$ adjacent to $u$) so that
$v^{X_u} =\Gamma (u)$ has size $\frac{p^R-1}{k}$. We therefore have \ben \label{Eqa}
\frac{p^R-1}{k} \geq \frac{b}{\ell} \cdot (p^{a}-1)^{\ell}. \een

On the other hand, $X_{u}$ contains $G_u=W=\brac{\omb^k}$, and all $W$-orbits in $V \setminus \{u
\}$ have length $\frac{p^R-1}{k}$. It follows that all orbits of $X_{u}$ in $V \setminus \{ u \}$
have length a multiple of $\frac{p^R-1}{k}$. Now there exists an orbit of $X_{u}$ in $V \setminus
\{ u \}$ of length $b(p^{a}-1)$ (the set of $b$-tuples with exactly one entry different from
$\gamma$), and so 
\ben \label{Eqb} b(p^{a}-1) \geq \frac{p^R-1}{k}. 
\een 
Combining inequalities
(\ref{Eqa}) and (\ref{Eqb}), we obtain 
\ben \label{Eqc} 
\ell \geq (p^{a}-1)^{\ell - 1}. 
\een 
Since
$p^a \geq 5$, the inequality (\ref{Eqc}) holds if and only if $\ell =1$, and hence $\frac{p^R-1}{k}
= b(p^{a}-1)$ as claimed. This implies that $k = \frac{p^R-1}{b(p^a-1)} =
\frac{a(p^R-1)}{R(p^a-1)}$, which is a contradiction. Thus $X$ is a primitive subgroup of
$\AGL(R,p)$. \qedd

\section{The case where $k \mid (p-1)$: Proof of Theorem~\ref{thm:CYC}~(4)} 
\label{sec:subcase1}

Let $\Gamma = \GPq =\Cay(V,S)$, where $q=p^R$ and $V$, $S$ are as defined in 
Definition~\ref{defn:gpaley}, and suppose that $k$ divides $p-1$. 
Let $A:=\Aut(\Gamma)$. Recall from Section~\ref{sec:intro} that $A$ contains 
$X:=T \rtimes \brac{W,\al}$ as an arc-transitive subgroup, where $W =
\brac{\omb^k}$. We will prove that $A=X$. 

If $k=2$, then $\Gamma = \GPP{p^R}{2}$ is a Paley graph 
and, see for example \cite{Peisert2001}, $A=X$ . Thus we may assume that 
$k \geq 3$. Then, since $p-1\geq k \geq 3$, we have $p\geq 5$.
Suppose that  $R=ab$ with $b >1$. Then   
$\frac{p^R-1}{p^a-1} = (p^{a})^{b-1} + p^{R-2a}+ \cdots + p^a +1 >p^{a(b-1)}
\geq p>k$. Hence, by Theorem~\ref{thm:CYC}~(1), $\Gamma$ is connected and so,
by Theorem~\ref{thm:add}, $W$ is irreducible on $V$.
Also if $k=\frac{p^R-1}{b(p^a-1)}$, then 
\[
b= \frac{p^R-1}{k(p^a-1)}>\frac{p^R-1}{p(p^a-1)} > p^{ab-a-1} \geq p^{a(b-2)}.
\]
It follows, since $p \geq 5$, that $b=2$, $a=1$ and $k =\frac{p+1}{2}$. 
However, this contradicts the assumption $ k \mid (p-1)$. Thus it follows
from Theorem~\ref{thm:CYC}~(3) that $A= T \rtimes A_0 \leq \AGL(R,p)$, 
where $\brac{W, \al} \leq A_0 \leq\GL(R,p)$. Note that $A_0$ preserves 
$S \subset V$, and hence $A_0$ does not contain $\SL(R,p)$.

We identify $V=\F{p^R}$ with an $R$-dimensional vector space over the prime 
field $\Fp$. Let $a$ be minimal such that $a \geq 1$, $a \mid R$ and $A_0$ 
preserves on $V$ the structure of an $a$-dimensional vector space over a 
field of order $q_0 =p^{R/a}$. Then $A_0 \leq \Gamma L(a,q_0)$ 
acting on $V = \F{q_0}^a$. Let $Z:=\brac{\omb^{(q-1)/(q_0-1)}}=Z(\GL(a,q_0)) 
\cong \Z_{q_0-1}$.

\bl \label{lem:pre} 
If $a\leq2$ then $A=X$.
\el

\bpf
Suppose first that $a=1$. Then $A_0\leq \Gamma L(1,p^R) =\brac{\omb, \al}$.
Since $\brac{\omb}$ is regular on $V^*$, and since $A_0$
leaves $S=\brac{\om^k}$ invariant, it follows that $A_0 \cap \brac{\omb}= 
\brac{\omb^k}$. Hence $A_0=W$ and $A=X$. 

Suppose now that $a=2$. Consider the
canonical homomorphism $\varphi: \GL(2,q_0) \mapp \PGL(2,q_0)$, and for 
$H\leq\GL(2,q_0)$ let $\oline{H}:=\varphi(H)=HZ/Z$. 
Now $\oline{A_0} \nsupseteq \PSL(2,q_0)$ since $A_0\nsupseteq \SL(2,q_0)$. 
Also $\oline{W} =  WZ/Z$, and $WZ =\brac{\omb^k, \omb^{q_0+1}}$ has order 
$q_0^2-1$ if $k$ is odd and $\frac{q_0^2-1}{2}$ if $k$ is even. Hence 
\[ 
\oline{\brac{W, \al}} \cong \left\{
\begin{array}{ll}
                D_{2(q_0+1)} &\mbox{if $k$ is odd}\\
                D_{q_0+1}     &\mbox{if $k$ is even.}
                 \end{array}
                            \right.\]
It follows from the classification of the subgroups of $\PGL(2,q_0)$ and 
$\PSL(2,q_0)$ (see \cite[p. 417]{Suzuki1}) that either $\oline{A_0} \leq 
\oline{\brac{\omb, \al}} \cong D_{2(q_0+1)}$ or $\oline{A_0} \in 
\{A_4, S_4, A_5 \}$.
In the former case, $A_0 \leq \brac{\omb, \al}=\AGaL(1,q)$ and $a=1$, which is
a contradiction. In the latter case,  since $\oline{A_0} \geq
\Z_{(q_0+1)/2}$ and $p \geq 5$, it follows that $q_0=p=5$ or 7. 
Moreover since $\oline{A_0} \geq
D_{p+1}$ and $\oline{A_0} \nsupseteq \PSL(2,p)$, it follows that 
$\oline{A_0} =S_4$ and in both cases
$\oline{A_0}$ is transitive on 1-spaces. Thus $S$ consists of, say, 
$s$ points from each 1-space and $|S| = (p+1)s$. Since $|S| =
|V^*|/k$, we have $k=|V^*|/((p+1)s) = (p-1)/s$. Also, since $\Gamma$ is 
an undirected Cayley graph, $S=\brac{\om^k}=-S$, and hence $S$ contains $-1$
and $s\geq2$. Thus $k\leq(p-1)/2$, and since 
$k \geq 3$, we have $p =7$, $k =3$ and $s=2$. This is impossible since 
in this case $W=\brac{\omb^3}\cong \Z_{16}$ projects to $\oline{W} \cong 
\Z_8$ and $\oline{A_0} =S_4$ has no such subgroup. \epf

From now on we will assume that $a \geq 3$, $k \geq 3$, and $p \geq 5$. Then,
by an old result of Zsigmondy~\cite{Zsig} (or see \cite{GPPS}), there is a 
prime divisor $r$ of $p^R-1$ such that $r$ does not divide $p^c-1$ for 
any $c<R$. Then $p$ has multiplicative order $R$ modulo $r$, and in particular
$R$ divides $r-1$. Thus $r=Rs+1$ for some $s\geq1$. Such a prime $r$ is called 
a \emph{primitive prime divisor} of $p^R-1$.

Let $A_1 := A_0 \ \cap \GL(a,q_0)$.  
Since $W \subseteq \GL(a,q_0)$, it follows that $W
\subseteq A_1$, so $r$ divides $|A_1|$ and $A_1$ is irreducible. 
By the minimality of $a$, $A_1$ is not contained in a proper `extension field
subgroup' of $\GL(a,q_0)$. Also, 
since $k \mid (p-1)$, the order of $W$ is divisible
by $\frac{p^R-1}{p-1}$, and it follows that $W$, and hence also $A_1$, cannot 
be realised over a proper subfield of $\F{q_0}$. 
By \cite[Main Theorem]{GPPS} (noting that the groups in \cite[Examples 2.2,
2.3, 2.4]{GPPS} do not have all of these properties), 
either 
\newcounter{tt}
\begin{list} {{\rm \textbf{(\Alph{tt})}}}{\usecounter{tt}}
\item $A_1$ belongs to one of the families of Examples 2.1 or 2.5 in
\cite{GPPS}, or
\item $A_1$ is \textit{nearly simple}, that is, $L \leq A_1/(A_1 \cap Z)
\leq \Aut(L)$ for some nonabelian simple group $L$, with $L$ as in one of 
Examples 2.6 -- 2.9 in \cite{GPPS}.
\end{list}

\bl \label{lem:CaseA} The group $A_1$ satisfies condition $\mathbf{(B)}$. 
\el 

\bpf Suppose that $A_1$ is a subgroup of $\GL(a,q_0)$ in 
\cite[Example 2.1]{GPPS}. Then $A_1$ is a 
classical group containing $Y=\SL(a,q_0)$, or (for $a$ even) 
${\rm Sp}(a,q_0)$ or $\Omega^\pm(a,q_0)$, or (if $aq_0$ is odd) 
$\Omega(a,q_0)$, or (if $q_0$ is a square) ${\rm SU}(a,q_0)$. 
Since $A_1$ is not transitive on $V^*$, $A_1$ cannot contain  $\SL(a,q_0)$ 
or ${\rm Sp}(a,q_0)$. Also $A_1$ contains the irreducible element
$\omb^{p-1}$ of order $\frac{p^r-1}{p-1}$, whereas (see \cite{Aron} or 
\cite{Huppert70}) for the remaining groups $Y$, 
$|\brac{\omb} \ \cap \ N_{\GL(a,q_0)}(Y)| \leq (q_0^{a/2}+1)(q_0-1)$. 

Next suppose that $A_1$ is as in \cite[Example 2.5]{GPPS}. Then $a =2^m$, 
$r=a+1$, and $A_1$ is contained in $Z
\circ (S \cdot M_0)$ where $S$ and $M_0$ are as listed in Table~\ref{tab:no2}.
Also, since $r\geq R+1$, it follows that $a=R$, and so $r=R+1=2^m+1 \geq 5$ and
$q_0=p$. Elements in $S$ and $M_0$ have orders at most 4 and $2^{2m}-1 =R^2-1$
respectively (see \cite[Proof of Lemma 2]{Aron}), and $S \cap Z \cong \Z_2$. 
Thus elements of $A_1$ have order at most $2(R^2-1)(p-1)$, and as $A_1$
contains $\omb^{p-1}$, of order $\frac{p^R-1}{p-1}$, we have 
$\frac{p^R-1}{p-1} \leq2(R^{2}-1)(p-1)$. Since $p \geq 5$ and $R\geq4$,
\[
5^{R-2} -1 \leq  p^{R-2} -1 <  \frac{p^R-1}{(p-1)^2} <  2(R^2-1).
\]
Since $R \geq4$, this implies that $(R,p) = (4,5)$. However, $r=R+1 =5$ does
not divide $5^4-1$, contradicting the definition of $r$. \epf

\begin{table}[ht]
\begin{tabular}{c c c}
\hline $S$ & \ $M_0$ & \ $p$ \\ \hline  $4 \circ 2^{1+2m} = \Z_4 \circ D_8 \circ \cdots \circ D_8$
&
\ ${\rm Sp}(2m,2)$ & \ $p \equiv 1$ (mod 4) \\
$2_{-}^{1+2m} =D_8 \circ \cdots \circ D_8 \circ Q_8$ & \ ${\rm O}^{-}(2m,2)$ &  \\
$2_{+}^{1+2m} =D_8 \circ \cdots \circ D_8$ & \ ${\rm O}^{+}(2m,2)$ &  \\
\hline
\end{tabular}
\vs \caption{} \label{tab:no2}
\end{table}

Thus case (\textbf{B}) holds, and we need to consider the possibilities
for $A_1$ from Examples 2.6 - 2.9 in \cite{GPPS} (see also Tables~\ref{tab:no4}
- \ref{tab:no7})  with order divisible by the primitive prime divisor $r$ of
$p^R-1$, where $R=ab, a \geq 3$, $q_0=p^{b}\geq 5^{b}$. 
Here, $L \leq A_1/(A_1 \cap Z) \leq
\Aut(L)$ for some non-abelian simple group $L$. Note that in applying the 
results of \cite{GPPS},
the dimension $a$ is equal to $d=e$ in \cite{GPPS} and $b$ is the parameter 
$a$ in \cite{GPPS}. Also, most of the examples in Tables~\ref{tab:no4} -
\ref{tab:no7} have additional conditions for $q_0$, $p$, $r$ or $a$. We will not mention them here,
unless they are necessary for our calculations.

\begin{table}[ht]
\begin{tabular}{|c|cccccc|}
\hline \multicolumn{7}{|c|}{Example 2.6 (b) of \cite{GPPS}} \\
\hline
 $n$& 7&6&5&5&7&7\\
$a$ &4 &4&2&4&3&6\\
$r$ & 5&5&5&5&7&7\\
$b$ &1&1&2&1&2&1\\
\hline
\end{tabular}
\vs \caption{} \label{tab:no4}
\end{table}

\begin{table}[ht]
\begin{tabular}{|c|cccccccc|}
\hline \multicolumn{9}{|c|}{Example 2.7 of \cite{GPPS}} \\
\hline
$L$&$M_{11}$&$M_{12}$&$M_{22}$&$M_{23}$&$J_{2}$&$J_{3}$&$Ru$&$Suz$\\  
$a$&10      &10      &10      &22      &6      &18     &28  &12   \\
$r$&11      &11      &11      &23      &7      &19     &29  &13   \\
\hline
\end{tabular}
\vs \caption{} \label{tab:no5}
\end{table}

\begin{table}[ht]
\begin{tabular}{|c|cccc|}
\hline \multicolumn{5}{|c|}{Example 2.9 of \cite[Table 7]{GPPS}} \\
\hline
$L$&$G_{2}(4)$&${\rm PSU}(4,2)$&${\rm PSU}(4,3)$&$\PSL(3,4)$\\ 
$a$&12        &4               &6              &6            \\
$r$&13        &5               &7              &7            \\
\hline
\end{tabular}
\vs \caption{} \label{tab:no6}
\end{table}

\begin{table}[ht]
\begin{tabular}{|c|ccccc|}
\hline \multicolumn{6}{|c|}{Example 2.9 of \cite[Table 8]{GPPS}} \\
\hline
$L$&$\PSL(n,s)$           &${\rm PSU}(n,s)$&    ${\rm PSp}(2n,s)$&$\PSL(2,s)$&
$\PSL(2,s)$\\
   &$n \geq 3$, $n$ prime&$n \geq 3$, $n$ prime& $n =2^c \geq 2$&$s \geq 7$&
$s \geq 7$\\
$a$&$\frac{s^n-1}{s-1} -1$&$\frac{s^n+1}{s+1}-1$& $\frac{1}{2}(s^n-1)$&$s$, 
$s-1$, or  $\frac{1}{2}(s-1) $&$\frac{1}{2}(s-1)$\\
$r$&$a+1$                 &$a+1$& $a+1$& $a+1$& $2a+1$\\ \hline
\end{tabular}
\vs \caption{} \label{tab:no7}
\end{table}

\bl \label{lem:CB} 
The group $A_1$ does not satisfy condition  $\mathbf{(B)}$. 
\el 

\bpf Recall that $r\geq R+1\geq a+1\geq4$ and $p\geq 5$. 
Also $\omb^{p-1}\in A_1$ has order $\frac{p^R-1}{p-1}$, a multiple of $r$.
Let $\max(A_1)$ denote the maximum order of an element of $A_1$ having order 
a multiple of $r$. An easy calculation shows that 
\begin{equation}\label{order}
\max(A_1)\geq \frac{p^R-1}{p-1}>\left\{\begin{array}{ll}
2(p^2-1)(R+1)&\mbox{if either $R\geq5, p\geq5$, or $R=4, p\geq11$}\\
(p^2-1)(R+1)&\mbox{if $R=4, p=7$}\\
2(p-1)(R+1)&\mbox{if $R\geq4, p\geq5$.}\\
			 \end{array}\right.
\end{equation}

\noindent
{\sc Case} \cite[Example 2.6]{GPPS} \quad Suppose first that 
$A_n \leq A_1 \leq S_n \times Z$ with $n=a+1$ or $a+2$, and  $r=a+1$, so 
$R=a, q_0=p$. Then any element 
of $S_n$ of order a multiple of $r$ is an $r$-cycle, and therefore $\omb^{p-1}$
has order at most $r(p-1)$. Hence $(R+1)(p-1)=r(p-1) \geq
\frac{p^R-1}{p-1}$, contradicting (\ref{order}). Thus $L=A_n$ with $n, a, r, 
b$ as in one of the columns of 
Table~\ref{tab:no4} (see  \cite[Tables 2, 3 and 4
]{GPPS}). In all cases $r=ab+1=R+1$ and $r \geq n-\delta$, where $\delta=1$ 
except for column 1 where $\delta=2$. Thus an element of $S_n$
of order a multiple of $r$ has order at most $\delta r$. 
Hence an element of $A_1$ 
of order a multiple of $r$ has order at most $ \delta r(q_0-1)=\delta
(p^b-1)(R+1)$. By (\ref{order}), column 3 of Table~\ref{tab:no4} holds but 
with $r=p=5$, contradicting the definition of $r$.

\noindent
{\sc Case} \cite[Example 2.7]{GPPS}, see Table~\ref{tab:no5}, which contains
the examples from \cite[Table 5]{GPPS} for which $r$ is a primitive prime 
divisor of $p^R-1$. In all cases $q=p$ and $r=R+1$. An  
element of $\Aut(L)$ of order a multiple of $r$ has order at most $2r=2(R+1)$
(see the \textsc{Atlas} \cite{Atlas}), so $\max(A_1)\leq 2(p-1)(R+1)$, 
contradicting (\ref{order}).

\noindent
{\sc Case} \cite[Example 2.8]{GPPS}. \ These examples are listed in 
\cite[Table 6]{GPPS} and the only ones for which $r$ is a primitive 
prime divisor of $p^R-1$ are $L =G_2(q_0)$ with $(a,p)=(6,2)$ and
$L=Sz(q_0)$ with $(a,p) = (4,2)$. However these are not examples for us since
$p \geq 5$.

\noindent
{\sc Case} \cite[Example 2.9]{GPPS}, see Tables~\ref{tab:no5} and 
\ref{tab:no7}, which contain the examples from \cite[Table 7 and 8]{GPPS} 
respectively for which $r$ is a primitive prime divisor of $p^R-1$ and $p
\geq5$.  We deal with Table~\ref{tab:no6} first. Here $R=a$ and $r=R+1$. An  
element of $\Aut(L)$ of order a multiple of $r$ has order at most $\delta r$,
where $\delta$ is 1, 2, 4 and 3 for the columns of Table~\ref{tab:no6} 
respectively, (see the \textsc{Atlas} \cite{Atlas}). Thus 
$\max(A_1)\leq \delta(p-1)(R+1)$, contradicting (\ref{order}) in all four 
cases.

Now we turn to the examples in Table~\ref{tab:no7}. In all cases $\gcd(s,p)=1$,
and we have $s=s_0^c$ for some prime $s_0\ne p$ and $c\geq1$.
Following the notation used in \cite{Aron}, we
let $m(K)$ denote the maximum of the orders of the elements of a finite group 
$K$. Then by (\ref{order}), $
m(A_1)\geq\max(A_1)\geq \frac{p^R-1}{p-1}\geq\frac{q_0^a-1}{q_0-1}>q_0^{a-1}$.
Moreover, in all cases, $m(A_1)\leq (q_0-1)m(\Aut(L))$, so
\begin{equation}\label{sharp}
m(\Aut(L))>q_0^{a-2}.
\end{equation}
In column 1 of Table~\ref{tab:no7}, $a\geq s^2+s\geq6$, and
$m(\Aut(L))\leq m(\GL(n,s)) \cdot 2c = 2c(s^n-1)=2c(s-1)(a+1)<s(s-1)(a+1) 
<a(a+1)$
so by (\ref{sharp}), $q_0^{a-2}<a(a+1)$, which is a contradiction since 
$a\geq6, q_0\geq5$.

In column 2 of Table~\ref{tab:no7}, $a\geq s^2-s\geq s\geq c$, and
$m(\Aut(L))\leq m(\GL(n,s^2)) \cdot 2c = 2c(s^{2n}-1)<2c(s^n+1)^2=
2c(s+1)^2(a+1)^2\leq 2a(a+1)^4$. 
By (\ref{sharp}), $q_0^{a-2} <  2a(a+1)^4$
and, since $q_0 \geq 5$, this implies that $a \leq 9$. Also, as $r =
\frac{s^n+1}{s+1} = a+1$ is prime, $n$ is odd and $3\leq a \leq 9$, 
it follows that $(r,a,s,n)=(7,6,3,3)$ and $L ={\rm PSU}(3,3)$. By the 
\textsc{Atlas} \cite{Atlas}, $m(\Aut({\rm PSU}(3,3))=12$ and we have a 
contradiction to (\ref{sharp}).

In column 3 of Table~\ref{tab:no7}, $a\geq s+1> c$, and
$m(\Aut(L))\leq m(\GL(2n,s)) \cdot 2c = 2c(s^{2n}-1)=2c\cdot 2a(2a+2)<
8a^2(a+1)$. By (\ref{sharp}), $q_0^{a-2} <  8a^2(a+1)$
and, since $q_0 \geq 5$, this implies that $a \leq 6$. Also, as $r =
\frac{s^n+1}{2} = a+1$ is prime, $n\geq2$ and $3\leq a \leq 6$, 
it follows that $(r,a,s,n)=(5,4,3,2)$ and $L ={\rm PSp}(4,3)$. By the 
\textsc{Atlas} \cite{Atlas}, $m(\Aut({\rm PSp}(4,3))=12$ and we have a 
contradiction to (\ref{sharp}).

In columns 4 and 5 of Table~\ref{tab:no7}, $L =\PSL(2,s)$ with $s \geq
7$, $s \neq p$ and $a \geq \frac{s-1}{2} \geq 3$. Then $m(\Aut(L)) =s+1$ 
(see \cite{Aron}) and hence (\ref{sharp}) yields $q_0^{a-2} < s+1 \leq 2a+2$.  
Since $q \geq5$ and $a \geq\frac{s-1}{2} \geq 3$, this implies that $a =3 = 
\frac{s-1}{2}$ and $q_0 \leq 7$. Thus $L =
\PSL(2,7)$ and $p \neq s =7$. So $q_0=p=5$ (since $p \geq 5$) and $a=R=3$. 
However the only primitive prime divisor 
of $p^R-1 = 5^3-1$ is 31 whereas by Table~\ref{tab:no7}, $r \leq s=7$.\epf

It follows from the discussion and results of this section that 
Theorem~\ref{thm:CYC}~(4) is proved.

\section{Acknowledgements} \label{sec:ack}

This paper forms part of the PhD thesis of the first author completed at the University of Western
Australia under the supervision of Associate Professor Cai Heng Li and the second author. It also
forms part of the Discovery Grant DP0449429 funded by the Australian Research Council. The authors are grateful to Mikhail Klin for pointing out connections between their work and cyclotomic association schemes.

\end{document}